\documentclass[12pt]{article}

\usepackage{amsmath}
\usepackage{amssymb}
\usepackage[mathscr]{eucal}
\usepackage{graphicx}
\usepackage{color}

\newtheorem{Theorem}{\sc Theorem}

\newtheorem{Proposition}[Theorem]{\sc Proposition}

\newcommand{\cQ}{\mbox{{${\cal Q}$}}}

\newcommand{\cP}{\mbox{{${\cal P}$}}}

\newcommand{\R}{{\if mm {\rm I}\mkern -3mu{\rm R}\else \leavevmode
\hbox{I}\kern -.17em\hbox{R} \fi}}

\newcommand{\cC}{\mbox{{${\cal C}$}}}
\newcommand{\cL}{\mbox{{${\cal L}$}}}
\newcommand{\cV}{\mbox{{${\cal V}$}}}

\newcommand{\wf}{\mbox{{$\widetilde{f}$}}}
\newcommand{\wu}{\mbox{{$\widetilde{u}$}}}

\newcommand{\wQ}{\mbox{{$\widetilde{{\cal Q}}$}}}
\newcommand{\wK}{\mbox{{$\widetilde{K}$}}}
\newcommand{\wA}{\mbox{{$\widetilde{A}$}}}
\newcommand{\wP}{\mbox{{$\widetilde{{\cal P}}$}}}
\newcommand{\wL}{\mbox{{$\widetilde{{\cal L}}$}}}
\newcommand{\wV}{\mbox{{$\widetilde{{\cal V}}$}}}

\newcommand{\cF}{\mbox{{${\cal F}$}}}

\newcommand{\bu}{\mbox{\boldmath{$u$}}}
\newcommand{\bv}{\mbox{\boldmath{$v$}}}

\newcommand{\bx}{\mbox{\boldmath{$x$}}}

\newcommand{\fb}{\mbox{\boldmath{$f$}}}

\newcommand{\bvarepsilon}{\mbox{\boldmath{$\varepsilon$}}}

\newcommand{\bnu}{\mbox{\boldmath{$\nu$}}}

\newcommand{\bzero}{\mbox{\boldmath{$0$}}}

\def\sqr#1#2{{
    \vcenter{
         \vbox{\hrule height.#2pt
               \hbox{\vrule width.#2pt height#1pt \kern#1pt
                     \vrule width.#2pt
               }
               \hrule height.#2pt
         }
    }
}}

\def\real{\mathbb{R}}

\def\lista#1
{{ \itemindent 0.0cm \labelsep .2cm \leftmargin 0.8cm \rightmargin
0.0cm \labelwidth 0.6cm \topsep 0.0mm
\parsep 0.0mm
\itemsep 0.0mm
\begin{list}{}
{ \setlength{\leftmargin}{.8cm} \setlength{\rightmargin}{0.0cm}
\setlength{\parsep}{0.0mm} \setlength{\topsep}{.0mm}
\setlength{\parskip}{.0cm} \setlength{\itemsep}{.0cm} }
{#1}\end{list}} }

\newcounter{theorem}

\begin{document}


\title{\bf  Convergence Results for Optimal Control Problems Governed by Elliptic Quasivariational Inequalities}

\vspace{22mm}
{\author{Mircea Sofonea$^{1}$\footnote{Corresponding author, E-mail : sofonea@univ-perp.fr}\,\, and\, Domingo A. Tarzia$^{2}$\\[6mm]
{\it \small $^1$ Laboratoire de Math\'ematiques et Physique}\\
{\it \small
	University of Perpignan Via Domitia}
\\{\it\small 52 Avenue Paul Alduy, 66860 Perpignan, France}		\\[6mm]		
{\it\small $^2$  Departamento de Matematica-CONICET}\\ {\it \small Universidad Austral}\\
		{\it \small Paraguay 1950, S2000FZF Rosario, Argentina}}

\date{}
\maketitle
\thispagestyle{empty}

\vskip 6mm

\noindent {\small{\bf Abstract.}
We consider an optimal control problem $\cQ$ governed by an elliptic quasivariational inequality with unilateral constraints. The existence of optimal pairs of the problem is a well known result, see \cite{SS}, for instance.
We associate to $\cQ$ a new optimal control  problem  $\wQ$, obtained by perturbing  the state inequality (including the
set of constraints and the nonlinear operator) and the cost functional, as well. Then, we provide sufficient conditions which guarantee the convergence of solutions of   Problem $\wQ$ to a solution of Problem $\cQ$. The proofs are based on  convergence results for elliptic quasivariational inequalities, obtained by using arguments of compactness, lower semicontinuity, monotonicity, penalty  and various estimates. Finally, we illustrate the use of the abstract convergence results in the study of optimal control associated with two boundary value problems.
The first one describes the equilibrium of an elastic body in frictional contact  with an obstacle, the so-called foundation. The process is static and the contact is modeled with normal compliance and unilateral constraint, associated to a version of Coulomb's law of dry friction.  The second one describes a stationary heat transfer problem with unilateral constraints. For the two problems
we prove existence, uniqueness and convergence results  together with the corresponding physical interpretation.

}

\vskip 6mm


\vskip 6mm


\vskip 15mm

\section{Introduction}\label{s1}
\setcounter{equation}0


\medskip\noindent
The study of optimal control problems is motivated by important applications in Physics, Mechanics, Automatics and Systems Theory. For instance,
the control of mathematical models which describe the contact of deformable bodies, as well as their optimal shape design,  is of considerable theoretical and applied interest in Civil Engineering, Automotive Industry and Mechanics of Structures. Moreover, the control of the temperature field in heat transfer  proccesses is important in various industrial settings like metal forming, among others.

Most of the models in Physics, Mechanics and Engineering Science are expressed in terms of
strongly  nonlinear boundary value problems with partial  differential equations which, in a weak formulation, lead to variational inequalities.
The theory of variational inequalities was developed  based on
arguments of monotonicity and convexity, including properties of the subdifferential
of a convex function.  Because
of their importance in engineering applications, a considerable effort has been put into their analysis, control
and numerical simulations and the literature in the field is extensive.
Basic references in the field are
\cite{BC, B, G, Kind-St,Li}, for instance.
Results in the study of optimal control for
variational  and variational-hemivariational inequalities   have been discussed in several works,  including \cite{Ba,BT,F,L,Mig,MP, NST}  and \cite{PK,S}, respectively.
Applications of variational inequalities in Mechanics could be found in
the books \cite{C,DL,EJK,HS,HHNL,KO, P}, for instance.  Reference on optimal control  for inequality problems arising in Mechanics and Physics include \cite{BT1,BT2,Ca,MM1,MM2,MM3,SMig,SS,SX1}.

In this paper we consider an optimal control problem for a general class of elliptic quasivariational inequalities.  The functional framework is the following:  $X$ and $Y$ are real Hilbert spaces endowed with the inner products $(\cdot,\cdot)_X$ and  $(\cdot,\cdot)_Y$, respectively, $K\subset X$, $A:X\to X$, $j:X\times X\to\R$, $f\in Y$ and $\pi:X\to Y$. Then, the inequality problem we consider is the following.

\medskip\medskip\noindent{\bf Problem}  ${\cal P}$. {\it Find $u$ such that}
\begin{equation}\label{1}u\in K,\qquad(Au,v-u)_X+j(u,v)-j(u,u) \ge(f,\pi v-\pi u)_Y \qquad\forall\,v\in K.
\end{equation}

\medskip We  associate to Problem $\cP$ the set of admissible pairs defined by
\begin{equation}\label{2}
{\cal V}_{ad} = \{\,(u, f)\in K\times Y \ \mbox{such that}\   (\ref{1})\  \mbox{holds}\,\}
\end{equation}
and we consider a cost functional ${\cal L}:X\times Y\to\mathbb{R}$. Here and below, $X\times Y$  represents the product of the  Hilbert spaces $X$ and $Y$, equipped with the canonical inner product.
Then, the optimal control problem we study in this paper is the following.

\medskip\medskip\noindent
{\bf Problem} ${\cal Q}$.  {\it Find $(u^*, f^*)\in {\cal V}_{ad}$ such that}
\begin{equation}\label{3}
{\cal L}(u^*,f^*)=\min_{(u,f)\in {\cal V}_{ad}} {\cal L}(u,f).
\end{equation}

Next, consider a set $\wK\subset X$, an operator $\wA:X\to X$ and an element $\wf\in Y$. With these data we construct the following perturbation of Problem  $\cP$.

\medskip\medskip\noindent{\bf Problem}  $\wP$. {\it Find $\wu$ such that}
\begin{equation}\label{1w}\wu\in \wK,\qquad(\wA\wu,v-\wu)_X+j(\wu,v)-j(\wu,\wu) \ge(\wf,\pi v-\pi \wu)_Y \qquad\forall\,v\in \wK.
\end{equation}

\medskip
We associate to Problem $\wP$ the set of of admissible pairs  given by
\begin{equation}\label{2w}
\wV_{ad} = \{\,(\wu, \wf)\in \wK\times Y \ \mbox{such that}\   (\ref{1w})\  \mbox{holds}\,\}
\end{equation}
and, for  a cost functional $\wL:X\times Y\to\mathbb{R}$, we construct the following perturbation of the optimal control problem $\cQ$.

\medskip\medskip\noindent
{\bf Problem} $\wQ$.  {\it Find $(\wu^*, \wf^*)\in \wV_{ad}$ such that}
\begin{equation}\label{3w}
\wL(\wu^*,\wf^*)=\min_{(\wu,\wf)\in\wV_{ad}}\wL(\wu,\wf).
\end{equation}

The unique  solvability of problems $\cP$ and $\wP$, on one hand, and the solvability of problems $\cQ$ and $\wQ$, on the other hand, follow from well known results obtained in the literature, under appropriate assumptions on the data. Here, we shall use the existence and uniqueness results in \cite{SS}, which will be resumed in the next section.

Now,
a brief comparation between problems $\cP$ and $\wP$ shows that Problem $\wP$ is obtained from Problem $\cP$ by replacing the set $K$ with the set $\wK$, the operator $A$ with the operator $\wA$ and the element $f$ with $\wf$. A similar remark can be made concerning the optimal problems $\cQ$ and $\wQ$, in which the set $\cV_{ad}$ was replaced by the set  $\wV_{ad}$ and the functional $\cL$ was replaced with $\wL$. Therefore, since problems $\wP$ and $\wQ$ represent perturbations of $\cP$ and $\cQ$, respectively, a natural question is to establish the link between the solutions of these problems.

In this paper we  provide a partial answer to the question above. Our aim is three folds.  The first one is to formulate sufficient assumptions on the data which guarantee the convergence of the solution $\wu$ of Problem $\wP$ to the solution $u$ of Problem $\cP$. Our result in this matter is Theorem \ref{t1} below, which represents the first novelty of this paper.  Our second aim is to  prove that, under appropriate conditions, the solutions of Problem $\wQ$ converge  to a solution of Problem $\cQ$. Our  result in this matter is Theorem \ref{t2}, which represent the second novelty of this work. Finally, our third aim is to illustrate the use of these abstract results in the study of two relevant examples. The first one arises from Contact Mechanics and the second one describe a heat transfer process.


The rest of this manuscript  is structured as follows. In Section \ref{s2}  we resume the existence and uniqueness results in \cite{SS} obtained in the study of problems $\cP$ and $\cQ$.
Then, in Section \ref{s3} we state and prove our main result concerning the link between the solutions of problems $\cP$ and $\wP$, Theorem \ref{t1}.
In Section \ref{s4} we state and prove our main result concerning the link between the solutions of problems $\cQ$ and $\wQ$, Theorem \ref{t2}. The proofs of the theorems are based on  arguments of compactness, lower semicontinuity, monotonicity, penalty  and various estimates.
In Section \ref{s5} we illustrate these abstract results in the study of a mathematical model which describes the frictional contact of an elastic material with a rigid-deformable foundation.
The process is  static and the contact is described with normal compliance and unilateral constraint, associated to a version of  Coulomb's law of dry friction. We apply the abstract result in Sections \ref{s3} and \ref{s4} in the study of this problem and provide the corresponding mechanical interpretations. We end this paper with
Section \ref{s6} in which we prove that Theorems \ref{t1} and \ref{t2} can be used to obtain a version of our previous convergence results obtained in \cite{BT1}, in the study of a
heat transfer model with unilateral constraints.

\section{Problem statement and preliminaries}\label{s2}
\setcounter{equation}0

In Sections \ref{s2}--\ref{s4} below we use the functional framework described in the Introduction and we denote by  $\|\cdot\|_X$,  $\|\cdot\|_Y$ the norms on the spaces $X$ and $Y$, respectively. All the limits, upper and lower limits below are considered as $n\to\infty$, even if we do not mention it explicitly. The symbols ``$\rightharpoonup$"  and ``$\to$"
denote the weak and the strong convergence in various spaces which will be specified, except in the case when these convergence take place in $\real$.

\medskip
In the study of Problem $\cP$ we consider the following assumptions.
\begin{eqnarray}
&&\label{K}
\quad\ \ K \ \mbox{\rm is a nonempty, closed, convex subset of} \ X.
\\[2.5mm]
&&\label{A}
\left\{ \begin{array}{l}
A\ {\rm is\ a\ strongly\ monotone\ Lipschitz\ continuous\ operator,\, i.e.,}\\[0mm]
{\rm there\ exists}\ m>0\ {\rm and}\ M>0\ {\rm such\ that}\quad\ \ \\[2mm]
\mbox{\rm (a) }\ (Au - Av,u -v)_{X} \geq m \|u -v\|^{2}_{X}
\quad \forall\,u,\,v\in X,\\[2mm]
\mbox{\rm (b) }\ \|Au-Av\|_X\le M\, {\|u-v\|_X}\quad\forall\,u,\,v\in
X.
\end{array}\right.\\[2.5mm]
&&\label{j}\left\{\begin{array}{l}  
\mbox{(a) } {\rm For\ all}\ \eta\in X,\ j(\eta,\cdot):X\to \mathbb
R \mbox{ is
	convex  and lower semicontinuous}\\
	 \quad\ \mbox{(l.s.c.), }\\[2mm]
\mbox{(b) There exists  }\alpha\ge 0 \mbox{ such that }\\
\quad\quad j(\eta_1,v_2)-j(\eta_1,v_1)+j(\eta_2,v_1)-j(\eta_2,v_2)\\
\qquad\quad\leq \alpha\,\|\eta_1-\eta_2\|_X\|v_1-v_2\|_X\quad
\forall \,\eta_1,\,\eta_2,\,v_1,\,v_2\in X. \end{array}\right.\\[2.5mm]
&&\label{small}
\qquad\quad\ m>\alpha.
\\ [3mm]
&&\qquad\quad\label{f}f\in Y.
\\[3mm]
&&\label{pi}
\left\{ \begin{array}{l} \pi\ {\rm is\ a\ linear\ continuous\ operator,\ i.e.,}\\[0mm]
\quad	{\rm there\ exists} \ c_0>0\ {\rm such\ that}\\[0mm]
\quad\quad \|\pi v\|_{Y} \leq c_0\,\|v\|_X\quad \forall\,v\in X.
\end{array}\right.
\end{eqnarray}

\medskip
We now recall  the following  existence and uniqueness result, proved in \cite{SS}.

\begin{Theorem}\label{t0}
	Assume that $(\ref{K})$--$(\ref{pi})$  hold.
	Then, the quasivariational inequality $(\ref{1})$  has a unique solution.
\end{Theorem}

In the study of Problem ${\cal Q}$ we assume that
\begin{equation}\label{l}
{\cal L}(u,f)=g(u)+h(f)\qquad\forall\, u\in X,\ f\in Y,
\end{equation}
where $g$ and $h$ are functions which satisfy the following conditions.
\begin{equation}\label{g}
\left\{ \begin{array}{l}
g:X\to\R\ {\rm is\ continuous,\  positive\ and\ bounded,\,  i.e.,}\hspace{19mm}\\[2mm]
\mbox{\rm (a) }\ v_n\to v\quad{\rm in}\quad X \quad\Longrightarrow\quad g(v_n)\to g(v).
\\[2mm]
\mbox{\rm (b) }\ g(v)\ge 0\quad\forall\, v\in X.
\\[2mm]
\mbox{\rm (c) }\ g\ \mbox{maps bounded sets in}\ X\ \mbox {into bounded sets in}\ \mathbb{R}.\end{array}\right.
\end{equation}
\begin{equation}\label{h}
\left\{ \begin{array}{l}
h:Y\to\R\ {\rm is\ weakly\ lower\ semicontinuous\ and\ coercive, \ i.e.,\ \ }\\[2mm]
\mbox{\rm (a) }\ f_n\rightharpoonup f\quad{\rm in}\quad Y \quad\Longrightarrow\quad \displaystyle\liminf h(f_n)\ge h(f).\hspace{22mm}
\\[2mm]
\mbox{\rm (b) }\  \|f_n\|_Y\to\infty\quad\quad\ \Longrightarrow\quad h(f_n)
\to\infty.
\end{array}\right.
\end{equation}

\begin{eqnarray}
&&\label{jb}
\left\{ \begin{array}{l}
	\mbox{\rm  There exist  }\beta,\,\gamma\ge 0\ \mbox{\rm such that }\\
	\quad j(\eta,v_1)-j(\eta,v_2)\leq (\beta+\gamma\|\eta\|_X)\,\|v_1-v_2\|_X
	\quad\forall \,\eta,\,v_1,\,v_2\in X.\end{array}\right.
\\[3mm]
&&\label{smal}\quad\ m>\gamma.\\ [3mm]
&&\hspace{0mm}\label{jc}\left\{ \begin{array}{l}
\mbox{\rm  For any sequences  }\{\eta_n\}\subset X,\
\{u_n\}\subset X\ \mbox{\rm such that }\\
\quad \eta_n\rightharpoonup\eta\in X, \ u_n\rightharpoonup u\in X\quad{\rm one\ has}\\
\qquad\displaystyle\limsup\,\big(j(\eta_n,v)-j(\eta_n,u_n)\big)\le j(\eta,v)-j(\eta,u)\quad\forall\,v\in X.\ \ \ \end{array}\right.
\\[3mm]
&&\hspace{0mm}\label{pic}
\left\{ \begin{array}{l}{\rm  For\ any\ sequence }\ \{v_n\}\subset X\ \mbox{\rm such that }\\
\mbox\quad v_n\rightharpoonup v\quad{\rm in}\quad X \quad{\rm one\ has}\quad \pi v_n\to\pi v\quad{\rm in}\quad Y.
\end{array}\right.
\end{eqnarray}

\bigskip
The following existence result was obtained in \cite{SS}.

\begin{Theorem}\label{t00}
	Assume that  $(\ref{K})$--$(\ref{small})$, $(\ref{pi})$--$(\ref{pic})$,
	Then, there exists at least one solution $(u^*, f^*)\in {\cal V}_{ad}$ of Problem ${\cal Q}$.
\end{Theorem}

\medskip

The proofs of Theorems \ref{t0} and \ref{t00} are based on  arguments of compactness, lower semicontinuity and monotonicity.
We shall use these theorems  in Sections \ref{s3} and \ref{s4} below, in the study of specific perturbed versions of problems $\cP$ and $\cQ$.

\section{A convergence result}\label{s3}
\setcounter{equation}0

In this section we  state and prove a  convergence result for the solution of Problem $\wP$, in the case when this problem has a specific structure. To this end,
we consider two sequences $\{\lambda_n\}\subset\real$, $\{f_n\}\subset Y$ and an operator $G:X\to X$. For each $n\in\mathbb{N}$
let $A_n:X\to X$ be the operator defined by
\begin{equation}\label{An}
A_nu=Au+\frac{1}{\lambda_n}\,Gu\qquad\forall\,u\in X,
\end{equation}
and denote by $\cP_n$  the following  version of Problem $\wP$, obtained with $\wA=A_n$ and $\wf=f_n$.

\medskip\noindent{\bf Problem}  ${\cal P}_n$. {\it Find $u_n$ such that}
\begin{eqnarray}\label{z0}
&&u_n\in\wK,\quad (A u_n, v - u_n)_X +
\frac{1}{\lambda_n} (Gu_n, v - u_n)_X +
j(u_n,v) - j(u_n,u_n)\label{2n}\\[2mm]
&&\qquad \qquad\qquad
\ge (f_n, \pi v - \pi u_n)_Y
\qquad\forall\, v \in \wK.\nonumber
\nonumber
\end{eqnarray}

\medskip
Note that  in the case when $\wK=X$, under appropriate assumptions on $G$,  Problem ${\cal P}_n$ represents a penalty problem of $\cP$.
Penalty methods have been widely used in the literature as an approximation tool to treat constraints in variational inequalities, as explained in \cite{G,KO,SofMat} and the references therein.

\medskip
To prove  the  unique solvability of Problem $\cP_n$ we use the following assumptions.
\begin{eqnarray}
&&\label{z1}
 \wK \ \mbox{\rm is a nonempty, closed, convex subset of} \ X.\\ [2mm]
&&\label{z2}G:X\to X \ \ \mbox{is a  monotone  Lipschitz continuous operator}.\\ [2mm]
&&\label{z3}\lambda_n>0\qquad\forall\,n\in\mathbb{N}. \\ [2mm]
&&\label{z4} f_n\in Y\qquad\forall\,n\in\ \mathbb{N}.
\end{eqnarray}

\medskip
We have the following  existence and uniqueness result.

\begin{Proposition}\label{p1}
	Assume  $(\ref{A})$--$(\ref{small})$, $(\ref{pi})$,  $(\ref{z1})$--$(\ref{z4})$. Then, for each $n\in\mathbb{N}$, there exists a unique solution $u_n\in X$
	to Problem~${\cal P}_n$.
\end{Proposition}

\noindent{\it Proof.}
Let  $n\in\mathbb{N}$. Assumptions (\ref{A}), (\ref{z2}), (\ref{z3})  imply that the operator
$A_n$ satisfies inequality (\ref{A})(a) with the same  constant $m$ as the operator $A$  and, moreover, it is Lipschitz continuous.  We conclude from above that the operator $A_n$ satisfies condition (\ref{A}).  Recall also assumptions  (\ref{z1}) and (\ref{z4}) on $\wK$ and $f_n$, respectively. These properties allows us to use Theorem \ref{t0} with $\wK$, $A_n$ and $f_n$ instead of $K$, $A$ and $f$, respectively. In this way we obtain the unique solvability of the inequality (\ref{z0}) which concludes the proof.
\hfill$\Box$

\medskip

To study the behavior of the solution of Problem ${\cal P}_n$ as $n\to\infty$ we consider the following additional hypotheses.
\begin{eqnarray}
&&\label{z5}\quad \ \lambda_n\to 0\quad {\rm as}\quad n\to\infty.\\ [2mm]
&&\label{z6}\quad \ f_n\rightharpoonup f\quad{\rm in}\quad  Y\quad  {\rm as}\quad n\to\infty.\\ [2mm]
&&\label{z7}\quad \ K\subset\wK. \\[2mm]
\label{z8}
&&\left\{\begin{array}{ll}
\mbox{\rm (a)}
\quad(Gu,v-u)_X\le 0\qquad\forall\, u\in\wK,\ v\in K\\ [3mm]
\mbox{\rm (b)} \quad u\in \wK,\quad (Gu,v-u)_X=0\quad\forall\,v\in K\ \ \Longrightarrow\ \  u\in K.
\end{array}\right.
\end{eqnarray}

\medskip

Note that, in the case when $\wK=X$, condition (\ref{z8}) is satisfied for any penalty operator of the set $K$, see Definition 23 in \cite{SMBOOK} for details.

\medskip
Our main result in this section is the following.

\begin{Theorem}\label{t1}
	Assume  $(\ref{K})$--$(\ref{pi})$, $(\ref{jb})$--$(\ref{pic})$, $(\ref{z1})$--$(\ref{z8})$ and, for each $n\in\mathbb{N}$, denote by $u_n$ the solution of Problem
	${\cal P}_n$. Then
	$u_n\to u$ in $X$, as $n\to \infty$, where  $u$
	is the solution of Problem ${\cal P}$.
\end{Theorem}


\noindent{\it Proof.}
The proof of Theorem \ref{t1} is carried out in several steps. 

\medskip
\noindent {\it
	{\rm i)} A first weak convergence result}. We claim that there is an element ${\widetilde{u}} \in \wK$ and
	a subsequence  of $\{ u_n \}$,
	still denoted by $\{ u_n \}$,
	such that $u_n \rightharpoonup {\widetilde{u}}$ in $X$, as $n\to \infty$.

\medskip
To prove the claim, we  establish the boundedness of the sequence $\{ u_n \}$ in $X$. Let $n\in\mathbb{N}$. We use assumption (\ref{z7})
and take $v=u$ in (\ref{z0}) to see that
\begin{eqnarray*}
(A u_n, u_n - u)_X\le
\frac{1}{\lambda_n} (Gu_n, u - u_n)_X
+j(u_n,u) - j(u_n,u_n)+
(f_n, \pi u_n - \pi u)_Y.
\end{eqnarray*}
Then,  using
the strong monotonicity of the operator $A$ we obtain that
\begin{eqnarray}
&&\label{5}
m \, \| u_n  - u \|_X^2 \le
(A u, u - u_n)_X+
\frac{1}{\lambda_n} (Gu_n, u - u_n)_X  \\ [2mm]
&&
\quad +
j(u_n,u) - j(u_n,u_n)+
(f_n, \pi u_n - \pi u)_Y. \nonumber
\end{eqnarray}
Next, assumption (\ref{z8})(a) implies that
\begin{equation}\label{6}
(Gu_n, u - u_n)_X\le0,
\end{equation}
and assumptions (\ref{j}), (\ref{jb})  yield
\begin{eqnarray}
&&\hspace{-8mm}\label{6n}
j(u_n,u) - j(u_n,u_n)\\ [2mm]
&&\hspace{-5mm}	=\big(j(u_n,u) - j(u_n,u_n)+ 	j(u,u_n) - j(u,u)\big)+\big(j(u,u) - j(u,u_n))\nonumber \\ [2mm]
&&\hspace{-2mm} \le \alpha\| u_n  - u \|_X^2+(\beta+\gamma\| u\|_X) \| u_n  - u \|_X.\nonumber
\end{eqnarray}
On the other hand, using (\ref{pi})
we find that
\begin{equation}\label{8}
(A u, u - u_n)_X+(f_n, \pi u_n - \pi u)_Y \le\big(\| Au\|_X+c_0\|f_n \|_{Y}\big)\| u_n - u \|_X.
\end{equation}
We now combine inequalities (\ref{5})--(\ref{8}) to see that
\begin{eqnarray}
&&\label{9}
m \, \| u_n  - u \|_X^2 \le (\| Au\|_X+c_0\|f_n \|_{Y}\big)\| u_n - u \|_X \\ [2mm]
&&
\quad\qquad +\alpha\| u_n - u \|_X^2 + (\beta + \gamma \| u \|_X) \| u_n - u \|_X.\nonumber
\end{eqnarray}

Note that by (\ref{z6}) we know that the sequence $\{f_n\}$ is bounded in $Y$. Therefore, using inequality (\ref{9}) and
the smallness assumption (\ref{small}),
we deduce that  there exists a constant $C > 0$
independent of $n$ such that $\| u_n-u\|_X \le C$. This implies that  the sequence $\{u_n\}$ is bounded  in $X$.
Thus, from the reflexivity of $X$, by passing to
a subsequence, if necessary,  we deduce that
\begin{equation}\label{10}
u_n \rightharpoonup {\widetilde{u}} \  \ {\rm in} \ \ X, \ \ \mbox{as}\ \ n\to \infty,
\end{equation}
with some ${\widetilde{u}} \in X$. Moreover, assumption (\ref{z1}) and the convergence (\ref{10}) implies that $\wu\in\wK$ and  completes the proof of the claim.

\medskip\noindent {\it
	{\rm ii)}  A property of the weak limit.}
	Next, we show that ${\widetilde{u}}$ is a solution to
	Problem ${\cal P}$.

\medskip
Let $v$ be a given element in $\wK$ and let $n\in\mathbb{N}$. We use  (\ref{z0}) to obtain that
\begin{eqnarray}
&&\label{11b}
\frac{1}{\lambda_n} (Gu_n, u_n - v)_X\le (Au_n, v-u_n)_X\\[2mm]
&&\qquad+j(u_n,v)-j(u_n,u_n)+(f_n, \pi u_n - \pi v)_Y.\nonumber
\end{eqnarray}
Then, by conditions (\ref{A}), (\ref{z6}), (\ref{jb}),  (\ref{pi}), using the boundedness of the sequence  $ \{u_n\}$, we deduce that each term  in the right hand side of inequality (\ref{11b}) is bounded. This implies that there  exists a constant
$D>0$  which does not depend on $n$ such that
\begin{equation*}
(Gu_n, u_n - v)_X\le \lambda_n D.
\end{equation*}
We now pass to the upper limit in this inequality and use the convergence (\ref{z5}) to deduce that
\begin{equation}\label{13}
\limsup\,(Gu_n, u_n -v)_X \le 0.
\end{equation}
Next, we take $v=\wu$ in (\ref{13}) and find that
\begin{equation}\label{13n}
\limsup\,(Gu_n, u_n -\wu)_X \le 0.
\end{equation}
Therefore, using assumption (\ref{z2}) and a standard    pseudomonotonicity argument (Proposition 1.23 in \cite{SofMat})  we obtain that
\begin{equation}\label{13nn}
\liminf\,(Gu_n, u_n -v)_X \ge (G\wu,\wu-v)_X\qquad\forall\, v\in X.
\end{equation}
We now combine inequalities (\ref{13nn})  and (\ref{13}) to find that
$(G {\widetilde{u}}, {\widetilde{u}} - v)_X\le 0$ for all $v\in\wK$.
Using now assumption (\ref{z8})(b) we deduce that ${\widetilde{u}} \in K$.

Consider now an  element $v\in K$. We use (\ref{z7}) and (\ref{z0})    to obtain
that
\begin{eqnarray*}
	&&(A u_n, u_n-v)_X\le
	\frac{1}{\lambda_n} (Gu_n, v - u_n)_X\\[2mm]
	&&\qquad
	+ j(u_n,v) - j(u_n,u_n)
	+ (f_n, \pi u_n - \pi v)_Y.
\end{eqnarray*}
Therefore, using
assumption (\ref{z8})(a)
we find that
\begin{equation}\label{17}
(A u_n, u_n - v)_X \le
j(u_n,v) - j(u_n,u_n)  + (f_n, \pi u_n - \pi v)_Y.
\end{equation}
Next, using (\ref{10}) and assumption (\ref{jc}) we have
\begin{equation}\label{18m}
\limsup\, \big(j(u_n,v) - j(u_n,u_n)\big)\le j({\widetilde{u}},v) - j({\widetilde{u}},{\widetilde{u}}).
\end{equation}
On the othe hand, assumption (\ref{z6}), (\ref{pic})
and the convergence (\ref{10}) yield
\begin{equation}\label{20}
(f_n, \pi u_n - v)_X \to (f, \pi\widetilde{u} - \pi v)_Y.
\end{equation}

We now use relations (\ref{17})--(\ref{20}) to see that
\begin{equation}
\label{21}
\limsup\,(A u_n, u_n - v)_X \le j({\widetilde{u}},v) - j({\widetilde{u}},{\widetilde{u}})
 + (f, \pi\widetilde{u} - \pi v)_X.
\end{equation}
Now, taking $v= {\widetilde{u}} \in K$ in (\ref{21})
we obtain that
\begin{equation}\label{24n}
\limsup\, (A u_n, u_n - {\widetilde{u}})_X \le 0.
\end{equation}
This inequality together with (\ref{10}) and the pseudomonotonicity
of $A$ implies that
\begin{equation}\label{24}
(A {\widetilde{u}}, {\widetilde{u}} - v)_X \le
\liminf\, (A u_n, u_n - v)_X
\ \  \ \ \forall \,  v \in X.
\end{equation}

\noindent
Combining now (\ref{24}) and (\ref{21}), we have
\begin{equation*}
(A {\widetilde{u}}, {\widetilde{u}} - v)_X\\ [2mm]
\le j(\widetilde{u},v) - j({\widetilde{u},\widetilde{u}}) +
 (f, \pi{\widetilde{u}} - \pi v)_Y
\end{equation*}

\noindent
for all $v \in K$. Hence, it follows that ${\widetilde{u}} \in K$ is a solution to Problem~${\cal P}$, as claimed.

\medskip\noindent {\it
	{\rm iii)}  A second weak convergence result}. We now prove the weak convergence of the whole sequence $\{u_n\}$.

\medskip
Since Problem~${\cal P}$ has a unique solution $u \in K$, we deduce from the previous step that ${\widetilde{u}} = u$. Moreover, a careful analysis of the proof in step {\rm ii)} reveals that every subsequence of $\{ u_n \}$
which converges weakly in $X$ has the  weak limit $u$. In addition, we recall that the sequence  $\{ u_n \}$ is bounded in $X$. Therefore, using  a standard argument we deduce that  the whole sequence $\{ u_n \}$ converges weakly in $X$ to $u$,  as $n\to \infty$.

\medskip\noindent {\it
	{\rm iv)} Strong convergence. } In the final step of the proof, we prove that
	$u_n \to u$ in $X$, as $n\to\infty$.

\medskip
We take
$v = {\widetilde{u}} \in K$ in (\ref{24}) and use (\ref{24n}) to obtain
\[
0 \le \liminf\, (A u_n, u_n - {\widetilde{u}})_X
\le
\limsup\,(A u_n, u_n - {\widetilde{u}})_X \le 0,
\]
which shows that $(A u_n, u_n - {\widetilde{u}})_X \to 0$,
as $n\to\infty$.
Therefore, using equality $\widetilde{u}=u$, the strong monotonicity of $A$
and the convergence $u_n \rightharpoonup u$  in $X$,
we have
\[
m_A \| u_n - u \|_X^2 \le
(A u_n - Au, u_n - u)_X =
(A u_n, u_n - u)_X - (A u, u_n - u)_X \to 0,
\]
as $n\to \infty$. Hence, it follows that $u_n \to u$ in $X$, which completes the proof.
\hfill$\Box$

\section{Convergence of optimal pairs }\label{s4}
\setcounter{equation}0

In this section we associate to Problem $\cP_n$ an optimal control problem for which we prove a convergence result. To this end
we keep the notation and assumptions in the previous section and we define the set of admissible pairs for  Problem $\cP_n$ by
\begin{equation}\label{2r}
{\cal V}_{ad}^n = \{\,(u_n, f_n)\in \wK\times Y \ \mbox{such that}\ (\ref{z0})\ \mbox{holds}\,\}.
\end{equation}
Then, the optimal control problem associated to Problem ${\cal P}_n$ is the following.

\medskip\noindent
{\bf Problem} ${\cal Q}_n$.  {\it Find $(u^*_n, f^*_n)\in {\cal V}_{ad}^n$ such that}
\begin{equation}\label{3r}
{\cal L}_n(u^*_n,f^*_n)=\min_{(u_n,f_n)\in {\cal V}_{ad}^n} {\cal L}_n(u_n,f_n).
\end{equation}

\medskip

In the study of Problem ${\cal Q}_n$ we assume that
\begin{equation}\label{ln}
{\cal L}_n(u,f)=g_n(u)+h_n(f)\qquad\forall\, u\in X,\ f\in Y,
\end{equation}
where $g_n$ and $h_n$ are functions which satisfy assumptions (\ref{g}) and (\ref{h}), respectively, for each $n\in\mathbb{N}$.
Note than when we use these assumptions for the functions $g_n$ and $h_n$ we refer to them as assumption (\ref{g})$_n$ and (\ref{h})$_n$, respectively.
Using Theorem \ref{t00}  we have the following existence result.

\begin{Proposition}\label{p2}
	Assume that  $(\ref{A})$--$(\ref{small})$, $(\ref{pi})$, $(\ref{ln})$, $(\ref{g})_n$, $(\ref{h})_n$,
	 $(\ref{jb})$--$(\ref{pic})$ and $(\ref{z1})$--$(\ref{z4})$ hold.
	Then, for each $n\in\mathbb{N}$, there exists
	at least one solution $(u^*_n, f^*_n)\in {\cal V}^n_{ad}$ of Problem ${\cal Q}_n$.
\end{Proposition}

To study the behavior of the sequence  of solutions of Problems ${\cal Q}_n$ as $n\to\infty$ we consider the following additional hypotheses.
\begin{eqnarray}
&&\label{q1} u_n\to u\quad {\rm in}\quad X \quad\Longrightarrow\quad  g_n(u_n)\to g(u).\\ [2mm]
&&\label{q2} f_n\rightharpoonup f\quad {\rm in}\quad Y \quad\Longrightarrow\quad  \liminf\,h_n(f_n)\ge h(f).\\ [2mm]
&&\label{q3n} {\color{red}\|f_n\|_Y\to\infty\quad\quad\ \, \Longrightarrow\quad h_n(f_n)\to\infty.}\\ [2mm]
&&\label{q3} h_n(f)\to h(f)\qquad\forall\, f\in Y.
\end{eqnarray}

Our main  result in this section is the following.

\begin{Theorem}\label{t2}
	Assume that  $(\ref{K})$--$(\ref{small})$,
	$(\ref{pi})$--$(\ref{pic})$,
	$(\ref{z1})$--$(\ref{z5})$, $(\ref{z7})$, $(\ref{z8})$, $(\ref{g})_n$ $(\ref{h})_n$, $(\ref{ln})$--$(\ref{q3})$ hold
 and, moreover, assume that $\{(u_{n}^{*},f_{n}^{*})\}$ is a sequence of solutions of Problem ${\cal Q}_n$.  Then, there exists a subsequence of the sequence $\{(u_{n}^{*}, f_{n}^{*})\}$, again denoted by $\{(u_{n}^{*}, f_{n}^{*})\}$, and an element   $(u^*,f^*)\in X\times Y$ such that
	\begin{eqnarray}
	&&\label{se1}
	f_n^* \rightharpoonup f^*\quad\mbox{\rm in}\quad Y\quad\mbox{\rm as}\quad n\to \infty, \\[2mm]
	&&\label{se2}
	u_n^*\rightarrow u^*\quad \mbox{\rm in}\quad X\quad\mbox{\rm as}\quad n\to \infty, \\[2mm]
	&&\label{se3}
	(u^*,f^*)\quad {\rm is\ a\ solution\ of\ Problem}\ {\cal Q}.
	\end{eqnarray}
\end{Theorem}

\medskip\noindent {\it Proof.} The proof is carried out in several steps, as follows.

\medskip\noindent
{\rm i)} {\it A boundedness result.}
	We claim that the sequence $\{f_n^*\}$ is bounded in $Y$.
	
	\medskip
	Arguing by contradiction, assume  that $\{f_n^*\}$ is not bounded in $Y$. Then, passing to a subsequence still denoted $\{f_n^*\}$, we have
	\begin{equation}\label{8n}
	\|f_n^*\|_Y\to +\infty\quad\text{as}\quad n\to +\infty.
	\end{equation}
	We  use equality (\ref{ln}) and assumption (\ref{g})$_n$(b) to see that
	\begin{equation*}\label{Ldef1n}
	{\cal L}_n(u_n^*, f_n^*) \geq h_n(f_n^*).
	\end{equation*}
	Therefore,   passing to the  limit as $n\to\infty$ in this inequality and
	using  (\ref{8n}) {\color{red}combined with assumption (\ref{q3n})} we deduce that
	\begin{equation}\label{20z}
	\lim\,{\cal L}_n(u_n^*, f_n^*)= +\infty.
	\end{equation}
	
	On the other hand, since $(u_{n}^{*},f_{n}^{*})$ represents a solution to Problem ${\cal Q}_n$, for each $n\in\mathbb{N}$ we have
	\begin{equation}\label{21z}
	{\cal L}_n(u^*_n,f^*_n)\le{\cal L}_n(u_n,f_n)
	\qquad\forall\,(u_n,f_n)\in {\cal V}_{ad}^n.
	\end{equation}
	We now  denote by ${u}_n^0$ the solution of Problem ${\cal P}_n$ for $f_n=f$. Then $(u_n^0,f)\in{\cal V}_{ad}^n$ and, therefore, (\ref{21z}) and (\ref{ln}) imply that
	\begin{equation}\label{21n}
	{\cal L}_n(u^*_n,f^*_n)\le g_n(u_n^0)+h_n(f).
	\end{equation}
Note that the convergences (\ref{z5}) and (\ref{z6}) allows us to apply Theorem \ref{t1} in order to see that
\begin{equation}\label{21x}
	u_n^0\to u\quad \mbox{\rm in}\quad X\quad\mbox{\rm as}\quad n\to \infty	
\end{equation}
where, recall, $u$ represents the solution of Problem $\cP$. Then, assumptions (\ref{q1}) and (\ref{q3}}) imply that
\begin{equation}\label{22}
g_n(u_n^0)+h_n(f)\to g(u)+h(f).
\end{equation}
	Relations (\ref{20z}), (\ref{21n} ) and (\ref{22}) lead to a contradiction, which concludes the claim.
	
\medskip\noindent {\rm ii)} {\it Two convergence results.}	In this step we prove the convergences (\ref{se1}) and (\ref{se2}).

\medskip
 First, since the sequence
	$\{f_n^*\}$ is bounded in $Y$ we can find
	a subsequence again denoted by $\{f_n^*\}$ and an element $f^*\in Y$ such that (\ref{se1}) holds.
	Next, we denote by $u^*$ the solution of Problem ${\cal P}$ for $f=f^*$. Then, we have
	\begin{equation}\label{se6}
	(u^*, f^*)\in {\cal V}_{ad}.
	\end{equation}
	Moreover, assumption (\ref{z5}), the convergence (\ref{se1}) and Theorem \ref{t1}  imply that (\ref{se2}) holds, too.

	\medskip\noindent {\rm iii)} {\it Optimality of the limit.}
	We now prove that $(u^*,f^*)$ is a solution to the optimal control problem ${\cal Q}$.
	
	\medskip
	We use the convergences (\ref{se1}), (\ref{se2}) and assumptions (\ref{q1}), (\ref{q2}), to see that
	\begin{equation*}
	\liminf \big(g_n(u_n^*)+h_n(f_n^*)\big)\ge g(u^*)+h(f^*)
	\end{equation*}
	and, therefore,
	the structure (\ref{ln}) and (\ref{l}) of the functionals ${\cal L}_n$ and ${\cal L}$ shows that
	\begin{equation}\label{se8}
	{\cal L}(u^*,f^*)\leq\liminf\,{\cal L}_n(u_n^*,f_n^*).
	\end{equation}
	
	Next, we fix a  solution $({u}^*_0,{f}^*_0)$ of Problem ${\cal Q}$ and,
	in addition, for each $n\in\mathbb{N}$  we denote by $\widetilde{u}_n^0$ the solution of Problem ${\cal P}_n$ for
	$f_n={f}^*_0$. It follows from here that $(\widetilde{u}_n^0,{f}^*_0)\in{\cal V}^n_{ad}$ and, by the optimality of the pair
	$(u_n^*,f_n^*)$, we have that
	\begin{equation*}
	{\cal L}_n(u_n^*,f_n^*)\leq{\cal L}_n(\widetilde{u}_n^0,{f}^*_0)\qquad\forall\, n\in\mathbb{N}.
	\end{equation*}
	We pass to the upper limit in this inequality to see that
	\begin{equation}\label{se9}
	\limsup\,{\cal L}_n(u_n^*, f_n^*)\leq \limsup\,{\cal L}_n(\widetilde{u}_n^0,{f}^*_0).
	\end{equation}

	Now, remember that ${u}^*_0$ is the
	solution of the inequality (\ref{1}) for $f=
	{f}^*_0$
	and $\widetilde{u}_n^0$  is the solution of the inequality (\ref{z0}) for $f_n=
	{f}^*_0$. Therefore, the convergence (\ref{z5}) and Theorem \ref{t1} imply that
	\begin{equation*}
	\widetilde{u}_n^0 \to {u}^*_0\quad\text{in}\quad X\quad\text{as}\quad n\to\infty
	\end{equation*}
	and, using  assumptions (\ref{q1}) and (\ref{q3}), we find that
	\begin{equation}\label{es}
	g_n(\widetilde{u}_n^0) \to g({u}^*_0),\quad h_n({f}^*_0)\to h({f}^*_0)\quad {\rm as}\quad n\to\infty.
	\end{equation}
	We now use (\ref{ln}),  (\ref{es}) and  (\ref{l}) to deduce that
	\begin{equation}\label{se10}
	\lim\,{\cal L}_n(\widetilde{u}_n^0,{f}^*_0)={\cal L}({u}^*_0,{f}^*_0).
	\end{equation}
	Therefore,  (\ref{se8}), (\ref{se9}) and (\ref{se10}) imply that
	\begin{equation}\label{se10n}
	{\cal L}(u^*,f^*)\leq {\cal L}({u}^*_0,{f}^*_0).
	\end{equation}
	
	On the other hand, since $({u}^*_0,{f}^*_0)$ is a solution of Problem ${\cal Q}$,  we have
	\begin{equation}\label{3p}
	{\cal L}({u}^*_0,{f}^*_0)=\min_{(u,f)\in {\cal V}_{ad}} {\cal L}(u,f).
	\end{equation} and, therefore,
	inclusion (\ref{se6})
	implies that
	\begin{equation}\label{se11}
	{\cal L}({u}^*_0,{f}^*_0)\le {\cal L}(u^*,f^*).
	\end{equation}
	We now combine the inequalities  (\ref{se10n}) and (\ref{se11}) to see that
	\begin{equation}\label{se16}
	{\cal L}(u^*, f^*)={\cal L}({u}^*_0,{f}^*_0).
	\end{equation}
	Finally, relations (\ref{se6}), (\ref{se16}) and (\ref{3p}) imply that (\ref{se3}) holds, which completes the proof of the Theorem.  \hfill$\Box$

\section{A frictional contact problem}\label{s5}
\setcounter{equation}0

The abstract results  in Sections \ref{s2}--\ref{s4} are useful
in the study of various mathematical models which describe the equilibrium of elastic bodies in frictional contact with a foundation. In this section
we provide an example of such model and, to this end, we need some  notations and preliminaries.

Let $d\in\{2,3\}$. We denote by $\mathbb{S}^d$ the space of second order symmetric tensors on $\mathbb{R}^d$ and use  the notation $``\cdot"$, $\|\cdot\|$, $\bzero$ for the inner product, the norm and the zero element of the spaces
$\mathbb{R}^d$ and $\mathbb{S}^d$, respectively.
Let $\Omega\subset\mathbb{R}^d$ be a domain with smooth boundary $\partial\Omega$ divided into three
measurable disjoint parts $\Gamma_1$, $\Gamma_2$ and $\Gamma_3$ such that ${ meas}\,(\Gamma_1)>0$.
A generic point in $\Omega\cup\Gamma$ will be denoted by $\bx=(x_i)$.
We use the
standard notation for Sobolev and Lebesgue spaces associated to
$\Omega$ and $\Gamma$. In particular, we use the spaces  $L^2(\Omega)^d$, $L^2(\Gamma_2)^d$,
$L^2(\Gamma_3)$  and $H^1(\Omega)^d$, endowed with their canonical inner products and associated norms.
Moreover, for an element $\bv\in H^1(\Omega)^d$ we still  write $\bv$ for the trace of
$\bv$ to $\Gamma$. In addition, we consider the
space
\begin{eqnarray*}
	&&V=\{\,\bv\in H^1(\Omega)^d\ :\  \bv =\bzero\ \ {\rm on\ \ }\Gamma_1\,\},
\end{eqnarray*}
which is a real Hilbert space
endowed with the canonical inner product
\begin{equation}
(\bu,\bv)_V= \int_{\Omega}
\bvarepsilon(\bu)\cdot\bvarepsilon(\bv)\,dx
\end{equation}
and the associated norm
$\|\cdot\|_V$. Here and below $\bvarepsilon$
represents the deformation operator, i.e.,
\[
\bvarepsilon(\bu)=(\varepsilon_{ij}(\bu)),\quad
\varepsilon_{ij}(\bu)=\frac{1}{2}\,(u_{i,j}+u_{j,i}),
\]
where an index that follows a comma denotes the
partial derivative with respect to the corresponding component of $\bx$, e.g.,\ $u_{i,j}=\frac{\partial u_i}{\partial j}$.
The completeness of the space $V$ follows from the assumption
${ meas}\,(\Gamma_1)>0$ which allows us to use Korn's inequality.
We denote by $\bzero_V$ the zero element of $V$ and we recall that, for an element $\bv\in V$,  the  normal and tangential components on $\Gamma$
are given by
$v_\nu=\bv\cdot\bnu$ and $\bv_\tau=\bv-v_\nu\bnu$, respectively.
We also  recall the trace inequality
\begin{equation}\label{trace}
\|\bv\|_{L^2(\Gamma)^d}\leq d_0\|\bv\|_{V}\qquad \forall\,
\bv\in V
\end{equation}
in which $d_0$ represents a positive constant.

For the inequality problem we consider in this section we use the data ${\cal F}$, $p$, $\fb_0$, $\fb_2$, $\mu$  and $k$ which satisfy
the following conditions.
\begin{eqnarray}
&&\left\{\begin{array}{ll}
{\rm (a)}\ {\cal F}\colon
\mathbb{S}^d\to \mathbb{S}^d. \\ [1mm]
{\rm (b)\  There\ exists}\ L_{\cal F}>0\ {\rm such\ that}\\
{}\qquad \|{\cal F}\bvarepsilon_1-{\cal F}\bvarepsilon_2\|
\le L_{\cal F} \|\bvarepsilon_1-\bvarepsilon_2\|\quad\mbox{for all} \ \ \bvarepsilon_1,\bvarepsilon_2
\in \mathbb{S}^d.
\\ [1mm]
{\rm (c)\ There\ exists}\ m_{\cal F}>0\ {\rm such\ that}\\
{}\qquad ({\cal F}\bvarepsilon_1-{\cal F}\bvarepsilon_2)
\cdot(\bvarepsilon_1-\bvarepsilon_2)\ge m_{\cal F}\,
\|\bvarepsilon_1-\bvarepsilon_2\|^2\quad \mbox{for all} \ \ \bvarepsilon_1,
\bvarepsilon_2 \in \mathbb{S}^d.\\ [1mm]
\end{array}\right.
\label{F}
\\[3mm]
&&\left\{\begin{array}{ll}
{\rm (a)}\ p:\mathbb{R}\to\mathbb{R}_+.\\ [1mm]
{\rm (b)\  There\ exists}\ L_p>0 {\rm\ such\
	that}\\
\qquad |p(r_1)-p(r_2)|\le L_p|r_1-r_2|\quad  \mbox{for all}\ \ r_1,\,r_2\in \mathbb{R}.\quad\ \\ [1mm]
{\rm (c)\ }
(p(r_1)-p(r_2))\,(r_1-r_2)\ge 0
\quad \mbox{for all}\ \
r_1,\,r_2\in \mathbb{R}.\\ [1mm]
{\rm (c)}\ p(r)=0 \quad {\rm iff}\quad r\le 0.
\end{array}\right.
\label{p}
\end{eqnarray}
\begin{eqnarray}\label{f0}
&&\fb_0\in L^2(\Omega)^d,\qquad\fb_2\in L^2(\Gamma_2)^d.\\ [2mm]
&&\label{mu}
\mu>0.\\ [2mm]
&&d_0^2\mu L_p<m_{\cal F}. \label{al}
\\ [2mm]
&&k>0. \label{k}
\end{eqnarray}
Moreover, we use 
$Y$ for the product space $L^2(\Omega)^d\times L^2(\Gamma_3)^d$ equipped with the canonical inner product,  and $K$ for the set defined by
\begin{equation}
\label{KK}K=\{\,\bv\in V\ :\ v_\nu \le k\ \  \hbox{a.e. on}\
\Gamma_3\,\}.
\end{equation}
Then, the inequality problem we consider  in this section is the following.

\medskip\medskip\noindent
\medskip\noindent{\bf Problem}  ${\cal P}^c$. {\it Find  $\bu$
	such that}
\begin{eqnarray}\label{51}
&&\bu\in K,\quad \int_{\Omega}{\cal F}\bvarepsilon(\bu)\cdot(\bvarepsilon(\bv)-\bvarepsilon(\bu))\,dx+\int_{\Gamma_3}p(u_\nu)(v_\nu-u_\nu)\,da\\[2mm]
&&\quad+\int_{\Gamma_3}\mu\,p(u_\nu)(\|{\bv}_\tau\|-\|{\bu}_\tau\|)\,da\ge
\int_{\Omega}\fb_0\cdot(\bv-\bu)\,dx
+\int_{\Gamma_2}\fb_2\cdot(\bv-\bu)\,da\quad\forall\,\bv\in K.\nonumber
\end{eqnarray}

Following the arguments in \cite{SofMat,SMBOOK}, it can be shown that Problem $\cP^c$ represents the variational formulation of a mathematical model that describes the equilibrium of an elastic body  $\Omega$ which is acted upon by  external forces,
is fixed on $\Gamma_1$, and is
in frictional contact on $\Gamma_3$. The contact takes place with a rigid foundation  covered by a layer of deformable material of thickness $k$.  In (\ref{51}) and below we shall refer to this foundation as foundation $F_k$.
Here ${\cal F}$ is the elasticity operator,   $\fb_0$ and $\fb_2$ denote the density of  applied body forces and tractions which act on the body and the surface $\Gamma_2$, respectively, $p$ is a given function which describes the reaction of the deformable material and $\mu$ represents the coefficient of friction.

\medskip
Next, we consider the constants $a_0$, $a_2$, $a_3$  and a function $\theta$ such that
\begin{eqnarray}
&&\label{52}a_0>0,\qquad a_2>0,\qquad a_3>0,\qquad \theta\in L^2(\Gamma_3).
\end{eqnarray}
We associate to Problem ${\cP}^c$ the set of admissible pairs ${\cal V}_{ad}^c$ and  the cost functional ${\cal L}$ given by
\begin{eqnarray}
&&\label{52b}
\hspace{-12mm}{\cal V}_{ad}^c = \{\,(\bu,\fb)\in K\times Y \ \mbox{such that}\  \fb=(\fb_0,\fb_2)\in Y\ \mbox{and}\ (\ref{51})\  \mbox{holds}\,\},\\[2mm]
&&\label{70}\hspace{-12mm}{\cal L}(\bu,\fb)=a_0\int_{\Omega}\|\fb_0\|^2\,dx+a_2\int_{\Gamma_2}\|\fb_2\|^2\,da+a_3\int_{\Gamma_3}|u_\nu-\theta|^2\,da
\end{eqnarray}
for all  $\bu\in V$, $\fb=(\fb_0,\fb_2)\in Y$.  Moreover, we consider the following optimal control problem.

\medskip\medskip\noindent
{\bf Problem} ${\cQ}^c$.  {\it Find $(\bu^*, \fb^*)\in {\cal V}^c_{ad}$ such that}
\begin{equation}\label{53}
{\cal L}(\bu^*,\fb^*)=\min_{(\bu,\fb)\in {\cal V}^c_{ad}} {\cal L}(\bu,\fb).
\end{equation}

Next, we consider a function $q$ and a constant $\widetilde{k}$ which satisfy the following conditions.
\begin{equation}
\left\{\begin{array}{ll}
	{\rm (a)}\ q:\mathbb{R}\to\mathbb{R}_+.\\ [1mm]
	{\rm (b)\ \ there\ exists}\ L_q>0 {\rm\ such\
		that}\\
	\qquad |q(r_1)-q(r_2)|\le L_q|r_1-r_2|\quad  \mbox{for all}\ \ r_1,\,r_2\in \mathbb{R}.\quad\ \\ [1mm]
	{\rm (c)\ }
	(q(r_1)-q(r_2))\,(r_1-r_2)\ge 0
	\quad \mbox{for all}\ \
	r_1,\,r_2\in \mathbb{R}.\\ [1mm]
	{\rm (d)}\ q(r)=0 \quad {\rm iff}\quad r\le 0.
\end{array}\right.
\label{q}
\end{equation}
\begin{equation}\label{kk}
\widetilde{k}\ge k>0.
\end{equation}
We introduce the set
\begin{eqnarray}
&&\label{54} \wK=\ \{\,\bv\in V\ :\ v_\nu\le \widetilde{k}\ \ {\rm on}\ \ \Gamma_3\,\}
\end{eqnarray}
and we assume that for each $n\in\mathbb{N}$ the functions $\fb_{0n}$, $\fb_{2n}$, $\theta_n$ and the constant $\lambda_n$ are given and satisfy the following conditions:
\begin{eqnarray}
&&\label{55}
\fb_{0n}\in L^2(\Omega)^d, \quad \fb_{2n}\in L^2(\Gamma_3)^d,\\ [2mm]
&&\label{55n}\lambda_n>0,\qquad \theta_n\in L^2(\Gamma_3).
\end{eqnarray}
Then, for each $n\in\mathbb{N}$, we consider the following perturbation of Problem $\cP^c$.

\medskip\medskip\noindent{\bf Problem}  ${\cal P}^c_n$. {\it Find $\bu_n$ such that}
\begin{eqnarray}
&&\label{56}\bu_n\in\wK,\quad \int_{\Omega}{\cal F}\bvarepsilon(\bu_n)
\cdot(\bvarepsilon(\bv)-\bvarepsilon(\bu_n))\,dx+
\int_{\Gamma_3}\,p(u_{n\nu})(v_\nu-u_{n\nu})\,da\\ [2mm]
&&\quad+
\frac{1}{\lambda_n}\int_{\Gamma_3}q(u_{n\nu}-k)
(v_\nu-u_{n\nu})\,da+
\mu\int_{\Gamma_3}p(u_{n\nu})(\|{\bv}_\tau\|-\|{\bu}_{n\tau}\|)\,da\nonumber\\[2mm]
&&\qquad\ge\int_{\Omega}\fb_{0n}\cdot(\bv-\bu)\,dx
+\int_{\Gamma_2}\fb_{2n}\cdot(\bv-\bu_n)\,da\qquad\forall\,\bv\in\wK.\nonumber
\end{eqnarray}

Following \cite{SofMat,SMBOOK}, Problem   ${\cal P}^c_n$ represents the variational formulation of the contact problem  with a foundation made of a rigid body covered by a layer of  deformable material of thickness $\widetilde{k}$. This layer is divided into two parts:
a first layer of thickness $\widetilde{k}-k>0$\, located on the top of the rigid body,
and a second layer of thickness $k$, located above.
Here, $\lambda_n$ is the deformability coefficient of the first layer and, therefore,  $\frac{1}{\lambda_n}$ represents its stiffness coefficient. In addition, $q$ is a given normal compliance function which describes the reaction of this first  layer.
We shall refer to this foundation as foundation $F_{\widetilde{k}}$.  A short comparation between the variational inequalities (\ref{51}) and (\ref{56})
reveals the fact that replacing the foundation $F_{{k}}$ with foundation $F_{\widetilde{k}}$ give rise to an extra term in the corresponding variational formulation, governed by the stiffness coefficient $\frac{1}{\lambda_n}$\,.

We  associate to Problem ${\cP}^c_n$ the set of admissible pairs ${\cal V}_{ad}^{cn}$ and  the cost function ${\cal L}_n$ given by
\begin{eqnarray}
&&\label{57}
\hspace{-12mm}{\cal V}_{ad}^{cn} = \{\,(\bu_n,\fb_n)\in\wK\times Y \ \mbox{such that}\  \fb=(\fb_{0n},\fb_{2n})\ \mbox{and}\ (\ref{56})\  \mbox{holds}\,\},\\[2mm]
&&\label{70n}\hspace{-12mm}{\cal L}_n(\bu_n,\fb_n)=a_0\int_{\Omega}\|\fb_{0n}\|^2\,dx+a_2\int_{\Gamma_2}\|\fb_{2n}\|^2\,da+a_3\int_{\Gamma_3}|u_{n\nu}-\theta_n|^2\,da
\end{eqnarray}
for all  $\bu_n\in V$, $\fb_n=(\fb_{0n},\fb_{2n})\in Y$.

\medskip

Our main result in this section, which represents a continuation of our previous results in \cite{SS}, is the following.

\begin{Theorem}\label{t5}
	Assume that $(\ref{F})$--$(\ref{k})$, $(\ref{52})$, $(\ref{q})$, $(\ref{kk})$, $(\ref{55})$ and  $(\ref{55n})$ hold. Then:
	
	\smallskip
	
	{\rm a)} Problem $\cP^c$  has a unique solution  and, for each $n\in\mathbb{N}$, Problem $\cP^c_n$ has a unique solution.
	Moreover, if
	\begin{equation}\label{co1c}
	\lambda_n\to 0,\quad \fb_{0n}\rightharpoonup \fb_{0}\ \ {\rm in}\ \ L^2(\Omega),\quad \fb_{2n}\rightharpoonup \fb_2\ \ {\rm in}\ \ L^2(\Gamma_3)\quad{\rm as}\ \ n\to \infty,
	\end{equation}
	the solution of Problem $\cP^c_n$ converges to the solution of Problem $\cP^c$, i.e.,
	\begin{equation}\label{co2c}
	\bu_n\to \bu\quad {\rm in}\ \ V\quad{\rm as}\quad n\to \infty.
	\end{equation}
	
	{\rm b)}	Problem $\cQ^c$  has at least one  solution and, for each $n\in\mathbb{N}$, Problem $\cQ^c_n$ has at least one solution.
	Moreover, if
	\begin{eqnarray}
	&&\label{co3c}
	\lambda_n\to 0, \qquad \theta_n\to\theta \ \ {\rm in}\ \  L^2(\Gamma_3)\qquad{\rm as}\ \ n\to \infty
	\end{eqnarray}
	and $\{(\bu_{n}^{*},\fb_{n}^{*})\}$ is a sequence of solutions of Problem ${\cal Q}^c_n$, there exists a subsequence of the sequence $\{(\bu_{n}^{*}, \fb_{n}^{*})\}$, again denoted by $\{(\bu_{n}^{*}, \fb_{n}^{*})\}$, and a solution   $(\bu^*,\fb^*)$ of Problem\ ${\cal Q}^c$, such that
	\begin{eqnarray}
	&&\label{se2x}
	\fb_n^* \rightharpoonup \fb^*\quad\mbox{\rm in}\quad Y,\quad \bu_n^*\rightarrow \bu^*\quad{\rm in}\quad V\quad
	\mbox{\rm as}\quad n\to \infty.
	\end{eqnarray}
	
\end{Theorem}

\medskip\noindent{\it Proof.} 
We start with some additional notation. First, we denote by  $\pi:V\to Y$ the operator
$\bv\mapsto (\iota\bv,\gamma_2\bv$) where  $\iota:V\to L^2(\Omega)^d$ is the canonic embedding and $\gamma_2:V\to L^2(\Gamma_2)^d$ is the restriction to the trace map to $\Gamma_2$.
Next, we
consider the operators $A:V\to V$, $G:V\to V$, the function $j:V\times V\to\R$ and the element $\fb\in Y$ defined as follows:
\begin{eqnarray}
&&
\label{8b1}(A\bu,\bv)_V =\int_{\Omega}\cF\bvarepsilon(\bu)\cdot\bvarepsilon(\bv)\,dx+\int_{\Gamma_3}\,p(u_{\nu})v_\nu\,da,\\ [2mm]
&&
\label{8b2}(G\bu,\bv)_V =\int_{\Gamma_3}q(u_\nu-k)v_\nu\,da,\\ [2mm]
&&j\colon V\times V \to \real,
\label{8b3}\quad
j(\bu,\bv)=\mu\int_{\Gamma_3}\, p(u_\nu)\|\bv_\tau\|\,da,\\[2mm]
&&\label{8b4}\fb=(\fb_0,\fb_2),
\end{eqnarray}
for all $\bu,\bv\in V$. Then it is easy to see that
\begin{equation}\label{e1c}
\left\{\begin{array}{l}
\mbox{$\bu$ is a solution of Problem $\cP^c$ if and only if}\\ [2mm]
\bu\in K, \quad (A\bu,\bv-\bu)_V+j(\bu,\bv)-j(\bu,\bu)\ge (\fb,\pi\bv-\pi\bu)_Y\quad \forall\, v\in K.
\end{array}\right.
\end{equation}
Moreover, for each $n\in\mathbb{N}$, \begin{equation}\label{e2c}
\left\{\begin{array}{l}
\mbox{$\bu_n$ is a solution of Problem $\cP^c_n$ if and only if\qquad}\\ [2mm]
\bu_n\in\wK, \quad (A\bu_n,\bv-\bu_n)_V+\frac{1}{\lambda_n}\,(G\bu_n,\bv-\bu_n)_V+j(\bu_n,\bv)-j(\bu_n,\bu_n)\\
\qquad\qquad\qquad\ge (\fb_n,\bv-\bu)_Y\quad \forall\, \bv\in\wK.
\end{array}\right.
\end{equation}

We now proceed with the proof of the two parts of the theorem.

\medskip\noindent

a) We use the abstract results in Sections \ref{s2} and \ref{s3} with $X=V$, $Y=L^2(\Omega)^d\times L^2(\Gamma_2)^d$, $K$ and $\wK$ defined by (\ref{KK}) and (\ref{54}), respectively, $A$ defined by (\ref{8b1}), $G$ defined by (\ref{8b2}), $j$ defined by (\ref{8b3}) and $\fb$ given by (\ref{8b4}). It is easy to see that in this case conditions $(\ref{K})$--$(\ref{pi})$, $(\ref{z1})$--$(\ref{z8})$ are satisfied.

\medskip
For instance, using assumption (\ref{F})  we see that
\begin{eqnarray*}
	(A\bu - A\bv,\bu -\bv)_{V} \geq m_{\cal F} \|\bu -\bv\|^{2}_{V},\qquad \|A\bu-A\bv\|_V\le (L_{\cal F}+d_0^2L_p)\, {\|\bu-\bv\|_V}
\end{eqnarray*}
for all $\bu,\, \bv\in V$. Therefore, condition (\ref{A}) holds with $m=m_{\cal F}$.
Condition (\ref{j})(a) is obviously satisfied and,
on the other hand, an elementary calculation based on the  definition (\ref{8b3}) and the trace inequality (\ref{trace}) shows that
\begin{eqnarray*}
	&&j(\bu_1,\bv_2)-j(\bu_1,\bv_1)+j(\bu_2,\bv_1)-j(\bu_2,\bv_2)\le d^2_0\mu L_p\,\|\bu_1-\bu_2\|_V\|\bv_1-\bv_2\|_V\quad
\end{eqnarray*}
for all $\bu_1,\,\bu_2,\,\bv_1,\,\bv_2\in V$. Therefore, condition (\ref{j})(b) holds with $\alpha=d_0^2\mu L_p$.
Next, condition $(\ref{jb})$ holds with $\beta=0$ and $\gamma=d_0^2\mu L_p$  and, using (\ref{al}) it follows that the smallness conditions (\ref{small})
and (\ref{smal}), too. We also note that conditions (\ref{jc}), and (\ref{pic}) arise from standard compactness arguments and, finally, condition  (\ref{z8}) is a direct consequence of the definitions  (\ref{8b2}),  (\ref{54}) and  (\ref{KK}), combined with the properties (\ref{q}) of the function $q$.

 Therefore, we are in a position to apply Theorem \ref{t0} and Proposition \ref{p1} in order to deduce the existence of a unique solution of the variational inequalities in (\ref{e1c}) and (\ref{e2c}), respectively.  Moreover, if (\ref{co1c}) holds, by Theorem \ref{t1} we deduce the convergence (\ref{co2c}).
These results combined with (\ref{e1c}) and (\ref{e2c}) allows us to conclude the proof of the first part of the theorem.

\medskip\noindent

b) Next, we use the abstract results in Sections \ref{s2} and \ref{s4} in the functional framework already described above, with the functionals $\cL$ and $\cL_n$ given by
(\ref{70}) and (\ref{70n}), respectively.
It is easy to see that in this case conditions
$(\ref{K})$--$(\ref{small})$,
$(\ref{pi})$--$(\ref{pic})$,
$(\ref{z1})$--$(\ref{z4})$, $(\ref{g})_n$ $(\ref{h})_n$, $(\ref{ln})$--$(\ref{q3})$ hold, with an appropriate choice of the functions $g$, $h$, $g_n$ and $h_n$.
Therefore, we are in a position to apply Proposition \ref{p2} in order to deduce the existence of a  solution of the optimal control problems in $\cQ^c$ and $\cQ^c_n$,
and Theorem \ref{t2}   in order to prove  the convergence (\ref{se2x}), as well.
\hfill$\Box$

\medskip
We now end this section with the following mechanical  interpretation of Theorem \ref{t5}.

\medskip
{\rm i)}  The convergence result (\ref{co2c}) shows that the solution of the frictional contact with foundation $F_k$ can be approximated by the solution of the frictional contact problem with foundation $F_{\widetilde{k}}$, with a large  stiffness coefficient of the first layer of the deformable material.
In other words, if this layer is almost rigid,  then the solution of the corresponding contact problem  is close to the solution
of the contact problem in which this layer is perfectly rigid.

\medskip
{\rm ii)} The mechanical interpretation of  the optimal control Problem $\cQ^c$ is the following:
given a contact  process governed by the variational inequality (\ref{51}) with the  data   ${\cal F}$, $p$, $k$ and $\mu$ which satisfy condition (\ref{F}), (\ref{p}), (\ref{mu}), (\ref{al}) and (\ref{k}),
we are looking for a couple of applied forces $(\fb_0,\fb_2)\in L^2(\Omega)^d\times L^2(\Gamma_2)^d$ such that
the normal displacement  of the solution  on the contact surface  is as close as possible to the ``desired" displacement $\theta$. Furthermore, this choice has to fulfill a minimum expenditure condition.
Theorem \ref{t5} guarantees the  existence of at least one optimal couple of applied forces $(\fb_0^*,\fb_2^*)$. A similar comment can be made on the optimal control Problem $\cQ^c_n$. Finally, the optimal solutions of the  contact problem associated to foundation $F_{\widetilde{k}}$ converge (in the sense given by Theorem \ref{t5} c)) to an optimal solution of the contact problem  associated foundation $F_k$, as the stiffness coefficient of the first deformable layer goes to infinity.

\section{A heat transfer boundary value problem}\label{s6}
\setcounter{equation}0

In this section we apply the abstract results in Sections \ref{s2}--\ref{s4} in the study of a mathematical model which describes a heat transfer phenomenon.  The problem we consider represents a version of the problem already considered in \cite{BT1} and, for this reason, we skip the details. Its classical formulation  is the following.

\medskip\noindent{\bf Problem $\cC^t$}. {\it
	Find a  temperature field $u:\Omega\to\R$ such that}
\begin{eqnarray}
&&\label{d1}u\ge 0,\qquad \Delta u+f\le 0,\qquad u(\Delta u+f)=0\qquad{\rm a.e.\ in\ }\Omega,\\ [2mm]
&&\label{d2}u=0\hspace{18mm}{\rm a.e.\ on\ }\Gamma_1,\\ [2mm]
&&\label{d3}u=b\hspace{18mm}{\rm a.e.\ on\ }\Gamma_2,\\ [2mm]
&&\label{d4}-\frac{\partial u}{\partial\nu}=q\hspace{11mm}{\rm a.e.\ on\ }\Gamma_3.
\end{eqnarray}

\medskip
 Here, as in Section \ref{s5},  $\Omega$ is a bounded domain in $\real^d$ ($d=1,2,3$ in applications) with smooth boundary $\partial\Omega=\Gamma_1\cup\Gamma_2\cup\Gamma_3$ and outer normal unit $\bnu$. We assume that $\Gamma_1$, $\Gamma_2$, $\Gamma_3$ are disjoint measurable sets and, moreover, $meas\,(\Gamma_1)>0$. In addition, in (\ref{d1})--(\ref{d4}) we do not mention the dependence of the different functions on the spatial variable $\bx\in\Omega\cup\partial\Omega$.
 The functions $f$, $b$ and $q$ are given and will be described below. Here we mention that  $f$ represents the internal energy, $b$ is  the prescribed temperature  field on $\Gamma_2$ and $q$ represents the heat flux prescribed on $\Gamma_3$. Moreover, $\frac{\partial u}{\partial\nu}$ denotes the normal derivative of $u$ on $\Gamma_3$.

\medskip
For the variational analysis of Problem $\cC^t$  we consider the space
\[V=\ \{\,v\in H^1(\Omega)\ :\ v=0\quad{\rm on}\quad\Gamma_1\}.\]
We denote in what follows by $(\cdot,\cdot)_V$ the inner product of the space $H^1(\Omega)$ restricted to $V$ and by $\|\cdot\|_V$ the associated norm.
Since $meas \,(\Gamma_1)>0$, it is well known that $(V,(\cdot,\cdot)_V)$ is a real Hilbert space.
Next, we assume that
\begin{eqnarray}
&&\label{dw}
f\in L^2(\Omega),\quad b\in L^2(\Gamma_2),\quad q\in L^2(\Gamma_3),\\ [2mm]
&&\label{dwz} \mbox{there exists $v_0\in V$ such that $v_0\ge 0$\ \ in\ \  $\Omega$\ and\  $v_0=b$\ \ on\ \ $\Gamma_2$}
\end{eqnarray}
and, finally, we introduce the set
\begin{eqnarray}
&&\label{d6} K=\ \{\,v\in V\ :\ v\ge 0\ \ {\rm in}\ \ \Omega,\quad v=b\ \ {\rm on}\ \ \Gamma_2\,\}.
\end{eqnarray}
Note that assumption (\ref{dwz}) represents a compatibility assumption on the data $b$ which guarantees that the set $K$ is not empty.
Then, it is easy to see that the variational formulation of problem $\cC^t$, obtained by standard arguments, is as follows.

\medskip\medskip\noindent{\bf Problem}  ${\cal P}^t$. {\it Find $u$ such that}
\begin{equation}\label{1d}u\in K,\quad\int_{\Omega}\nabla u\cdot(\nabla v-\nabla u)\,dx+\int_{\Gamma_3}q(v-u)\,da\ge\int_\Omega f(v-u) \,dx \quad\forall\,v\in K.
\end{equation}

\medskip
We now introduce the set of admissible pairs for inequality (\ref{1d}) defined by
\begin{equation}\label{2d}
{\cal V}_{ad}^t = \{\,(u, f)\in K\times L^2(\Omega) \ \mbox{such that}\   (\ref{1d})\  \mbox{holds}\,\}.
\end{equation}
Moreover, we consider two constants $\omega$, $\delta$  and a function $\phi$ such that
\begin{eqnarray}
&&\label{d15}\omega>0,\qquad  \delta>0,\qquad \phi\in L^2(\Omega)
\end{eqnarray}
and, with these data, we associate to Problem $\cP^t$ the following optimal control problem.

\medskip\medskip\noindent
{\bf Problem} ${\cal Q}^t$.  {\it Find $(u^*, f^*)\in {\cal V}_{ad}^t$ such that}
\begin{equation}\label{3d}
\omega\int_{\Omega}(u^*-\phi)^2\,dx+\delta\int_{\Omega}{(f^*)}^2\,dx=\min_{(u,f)\in {\cal V}_{ad}^t}\Big\{\omega\int_{\Omega}(u-\phi)^2\,dx+\delta\int_{\Omega}f^2\,dx\Big\}.
\end{equation}

Next, we introduce the set
\begin{eqnarray}
&&\label{d7} \wK=\ \{\,v\in V\ :\ v\ge 0\ \ {\rm in}\ \ \Omega\,\}
\end{eqnarray}
and we assume that for each $n\in\mathbb{N}$ the functions $f_n$, $\phi_n$ and the constants $\lambda_n$, $\omega_n$, $\delta_n$, are given and satisfy the following conditions:
\begin{eqnarray}
&&\label{dwx}
f_n\in L^2(\Omega),\\ [2mm]
&&\label{d15x}\lambda_n>0,\qquad \omega_n>0,\qquad  \delta_n>0,\qquad \phi_n\in L^2(\Omega).
\end{eqnarray}
Then, for each $n\in\mathbb{N}$, we consider the following perturbation of Problem $\cP^t$.

\medskip\medskip\noindent{\bf Problem}  ${\cal P}^t_n$. {\it Find $u_n$ such that}
\begin{eqnarray}
&&\label{1dx}u_n\in\wK,\qquad\int_{\Omega}\nabla u_n\cdot(\nabla v-\nabla u_n)\,dx+\int_{\Gamma_3}q(v-u_n)\,da\\[2mm]
&&\qquad+\frac{1}{\lambda_n}
\int_{\Gamma_2}(u_n-b)(v-u_n)\,da \ge\int_\Omega f_n(v-u_n) \,dx \qquad\forall\,v\in \wK.\nonumber
\end{eqnarray}

Using standard arguments it is easy to see that  Problem   ${\cal P}^t_n$ represents the variational formulation of the following boundary value problem.

\medskip\noindent{\bf Problem ${\cC}^t_n$}. {\it
	Find a temperature field $u_n:\Omega\to\R$ such that}
\begin{eqnarray}
&&\label{d1n}u_n\ge 0,\qquad \Delta u_n+f_n\le 0,\qquad u_n(\Delta u_n+f_n)=0\qquad{\rm a.e.\ in\ }\Omega,\\ [2mm]
&&\label{d2n} u_n=0\hspace{35mm}{\rm a.e.\ on\ }\Gamma_1,\\ [2mm]
&&\label{d3n}-\frac{\partial u_n}{\partial\nu}=\frac{1}{\lambda_n}(u_n-b)\hspace{11mm}{\rm a.e.\ on\ }\Gamma_2,\\ [2mm]
&&\label{d4n}-\frac{\partial u_n}{\partial\nu}=q\hspace{29mm}{\rm a.e.\ on\ }\Gamma_3.
\end{eqnarray}

Note that Problem ${\cC}^t_n$ is obtained from Problem $\cC^t$ by replacing the Dirichlet boundary condition (\ref{d3}) with the Neumann boundary condition (\ref{d3n}) and  prescribing the internal energy $f_n$ in $\Omega$, instead of the internal energy $f$. Here $\lambda_n$ is a positive parameter, and its inverse $h_n=\frac{1}{\lambda_n}$ represents the heat transfer coefficient on the boundary $\Gamma_2$. In contrast to Problem  $\cP^t$ (in which the temperature is prescribed on $\Gamma_2$), in Problem $\cP^t_n$ this condition is replaced by a condition on the flux of the temperature, governed by a positive  heat transfer coefficient.

\medskip
The set of admissible pairs for inequality (\ref{1dx}) is defined by
\begin{equation}\label{2dx}
{\cal V}_{ad}^{tn} = \{\,(u_n, f_n)\in \wK\times L^2(\Omega) \ \mbox{such that}\   (\ref{1dx})\  \mbox{holds}\,\}
\end{equation}
and, moreover, the associated optimal control problem is the following.

\medskip\medskip\noindent
{\bf Problem} ${\cal Q}^t_n$.  {\it Find $(u^*_n, f^*_n)\in {\cal V}_{ad}^{tn}$ such that}
\begin{eqnarray}\label{3dx}
&&\omega_n\int_{\Omega}(u^*_n-\phi_n)^2\,dx+\delta_n\int_{\Omega}{(f^*_n)}^2\,dx\\ [3mm]
&&\qquad=\min_{(u,f)\in {\cal V}_{ad}^{tn}}\Big\{\omega_n\int_{\Omega}(u-\phi_n)^2\,dx+\delta_n\int_{\Omega}f^2\,dx\Big\}.\nonumber
\end{eqnarray}

\medskip  Our main result in this section is the following.

\begin{Theorem}\label{t7}
Assume that $(\ref{dw})$--$(\ref{dwz})$, $(\ref{d15})$, $(\ref{dwx})$   and  $(\ref{d15x})$ hold. Then:

\smallskip

{\rm a)} Problem $\cP^t$  has a unique solution  and, for each $n\in\mathbb{N}$, Problem $\cP^t_n$ has a unique solution. Moreover, if
\begin{equation}\label{co1}
\lambda_n\to 0\quad{\rm and}\quad f_n\rightharpoonup f\ \ {\rm in}\ \ L^2(\Omega)\quad{\rm as}\ \ n\to \infty,
\end{equation}
the solution of Problem $\cP^t_n$ converges to the solution of Problem $\cP^t$, i.e.,
\begin{equation}\label{co2}
u_n\to u\quad {\rm in}\ \ V\qquad\quad{\rm as}\quad n\to \infty.
\end{equation}

{\rm b)}	Problem $\cQ^t$  has at least one  solution and, for each $n\in\mathbb{N}$, Problem $\cQ^t_n$ has at least one solution. Moreover, the solution of Problem $\cQ^t$ is unique if $\phi=0_{L^2(\Omega)}$ and, for each $n\in \mathbb{N}$, the solution of Problem $\cQ^{tn}$ is unique, if $\phi_n=0_{L^2(\Omega)}$.

\smallskip
{\rm c)} Assume that
\begin{eqnarray}
&&\label{co3}
\lambda_n\to 0,\quad\omega_n\to\omega,
\quad\delta_n\to\delta,  \quad \phi_n\to\phi \ \ {\rm in}\ \  L^2(\Omega)\quad{\rm as}\ \ n\to \infty
\end{eqnarray}
and let $\{(u_{n}^{*},f_{n}^{*})\}$ be a sequence of solutions of Problem ${\cal Q}^t_n$.  Then, there exists a subsequence of the sequence $\{(u_{n}^{*}, f_{n}^{*})\}$, again denoted by $\{(u_{n}^{*}, f_{n}^{*})\}$, and a solution   $(u^*,f^*)$ of Problem\ ${\cal Q}^t$, such that
\begin{eqnarray}
&&\label{se1x}
f_n^* \rightharpoonup f^*\quad\mbox{\rm in}\quad L^2(\Omega),\quad u_n^*\rightarrow u^*\quad{\rm in}\quad V\quad
\mbox{\rm as}\quad n\to \infty.
\end{eqnarray}
Moreover, if $\phi=0_{L^2(\Omega)}$, then the whole sequence $\{(u_{n}^{*}, f_{n}^{*})\}$ satisfies $(\ref{se1x})$ where $(u^*, f^*)$ represents the unique solution of Problem $\cQ^t$.
\end{Theorem}

\noindent{\it Proof.} We start by introducing some notation which allow us to write the problems in an equivalent form. To this end, we denote by  $\pi:V\to L^2(\Omega)$ the canonical inclusion of $V$ in $L^2(\Omega)$. Moreover, we
consider the operators $A:V\to V$, $G:V\to V$ defined as follows:
\begin{eqnarray}
&&\label{d8} \ (Au,v)_V=\int_\Omega \nabla u\cdot\nabla v\,dx+\int_{\Gamma_3}qv\,da\qquad\forall\,u, v\in V,\\ [2mm]
&&\label{d9} \ (Gu,v)_V=\int_{\Gamma_2}(u-b)v\,da\qquad\forall\,u, v\in V.
\end{eqnarray}
\medskip
Then, it is easy to see that
\begin{equation}\label{e1}
\left\{\begin{array}{l}
\mbox{$u$ is a solution of Problem $\cP^t$ if and only if}\\ [2mm]
u\in K, \quad (Au,v-u)_V\ge (f,v-u)_{L^2(\Omega)}\quad \forall\, v\in K.
\end{array}\right.
\end{equation}
Moreover, for each $n\in\mathbb{N}$, \begin{equation}\label{e2}
\left\{\begin{array}{l}
\mbox{$u_n$ is a solution of Problem $\cP^t_n$ if and only if\qquad}\\ [2mm]
u_n\in\wK, \quad (Au,v-u)_V+\frac{1}{\lambda_n}(Gu_n,v-u_n)_V\\
\qquad\qquad\qquad\ge (f_n,v-u_n)_{L^2(\Omega)}\quad \forall\, v\in\wK.
\end{array}\right.
\end{equation}

Next, denote by $\cL:V\times L^2(\Omega)\to\R$ and $\cL_n:V\times L^2(\Omega)\to\R$ the cost functionals given by
\begin{eqnarray}
&&\label{d9n}\cL(u,f)=\omega\|u-\phi\|^2_
{L^2(\Omega)}+\delta\|f\|^2_
{L^2(\Omega)}\,,
\\ [2mm]
&&\label{d10}\cL_n(u,f)=\omega_n\|u-\phi_n\|^2_
{L^2(\Omega)}+\delta_n\|f\|^2_
{L^2(\Omega)}
\end{eqnarray}	
for all	$(u,f)\in V\times L^2(\Omega)$. Then, it is easy to see that
\begin{equation}\label{e3}
\left\{\begin{array}{l}
\mbox{$(u^*,f^*)$ is a solution of Problem $\cQ^t$ if and only if\qquad}\\ [2mm]
(u^*,f^*)\in\cV_{ad}^t\quad{\rm and}\quad \cL(u^*,f^*)=\displaystyle\min_{(u^*f^*)\in\cV_{ad}^t}\cL(u,f)
\end{array}\right.
\end{equation}
Moreover, for each $n\in\mathbb{N}$,
\begin{equation}\label{e4}
\left\{\begin{array}{l}
\mbox{$(u^*_n,f^*_n)$ is a solution of Problem $\cQ^t$ if and only if}\\ [2mm]
(u^*_n,f^*_n)\in\cV_{ad}^{tn}\quad{\rm and}\quad \cL_n(u^*_n,f^*_n)=\displaystyle\min_{(u^*f^*)\in\cV_{ad}^{tn}}\cL_n(u,f)
\end{array}\right.
\end{equation}

We now proceed with the proof of the two parts of the theorem.

 \medskip
 a) We use the abstract results in Sections \ref{s2} and \ref{s3} with $X=V$, $Y=L^2(\Omega)$, $K$ and $\wK$ defined by (\ref{d6}) and (\ref{d7}), respectively, $A$ defined by (\ref{d8}), $G$ defined by (\ref{d9}),
 and $j\equiv 0$. It is easy to see that in this case conditions $(\ref{K})$--$(\ref{pi})$, $(\ref{z1})$--$(\ref{z8})$ are satisfied. Therefore, we are in a position to apply Theorem \ref{t0} and Proposition \ref{p1} in order to deduce the existence of a unique solution of the variational inequalities in (\ref{e1}) and (\ref{e2}), respectively.  Moreover, by Theorem \ref{t1} we deduce the convergence (\ref{co2}).
 These results combined with (\ref{e1}) and (\ref{e2}) allows us to conclude the proof of the statement a) in Theorem \ref{t7}.

 \medskip
 b)  We use the abstract results in Sections \ref{s2} and \ref{s4} in the functional framework described above,  with the functionals $\cL$ and $\cL_n$ given by
 (\ref{d9n}) and (\ref{d10}), respectively.
  It is easy to see that in this case conditions
 $(\ref{K})$--$(\ref{small})$,
 $(\ref{pi})$--$(\ref{pic})$,
 $(\ref{z1})$--$(\ref{z4})$, $(\ref{g})_n$ $(\ref{h})_n$, $(\ref{ln})$ and $(\ref{q1})$--$(\ref{q2})$ hold, with an appropriate choice of the functions $g$, $h$, $g_n$ and $h_n$.
  Therefore, we are in a position to apply Theorem \ref{t00} and Proposition \ref{p2} in order to deduce the existence of a  solution of the optimal control problems in (\ref{e3}) and (\ref{e4}), respectively.

  The uniqueness of the solution of Problem $\cQ^t$ in the case $\phi=0_{L^2(\Omega)}$ follows from a strict convexity argument. Indeed, for any $f\in L^2(\Omega)$ let  $u(f)$ denote the solution of the variational inequality in (\ref{e1}). Then it  was proved in \cite{BT1} that the functional
 \[f\mapsto\cL(u(f),f)=\omega\|u(f)\|^2_
 {L^2(\Omega)}+\delta\|f\|^2_
 {L^2(\Omega)}\] is strictly convex and, therefore, the optimal control problem in (\ref{e3}) has a unique solution. The uniqueness of the solution of Problem $\cQ^t_n$ in the case $\phi_n=0_{L^2(\Omega)}$ follows from the same argument. These results combined with the equivalence results (\ref{e3}) and (\ref{e4}) allows us to conclude the proof of the statement b)  in Theorem \ref{t7}.

 \medskip
 c) The convergence (\ref{se1x}) is a direct consequence of  Theorem \ref{t2}. The convergence $(\ref{se1x})$ of  the whole sequence $\{(u_{n}^{*}, f_{n}^{*})\}$ in the case $\phi=0_{L^2(\Omega)}$ follows from
 a standard argument, since in this case Problem $\cQ^t$ has a unique solution.
 \hfill$\Box$

\medskip
We end this section with the following physical interpretation of Theorem \ref{t7}.

\medskip
{\rm i)} First, the solutions of Problems $\cP^t$ and $\cP^t_n$ represent weak solutions of the the heat transfer problems $\cC^t$ and $\cC^t_n$, respectively. Therefore,
Theorem \ref{t7} provides the unique weak solvability of these problems. Moreover,
the weak solution of the problem with prescribed temperature on $\Gamma_2$ can be approximated by the solution of the problem with heat transfer on $\Gamma_2$, for a large heat transfer coefficient, as shown in \cite{DT}.

\medskip
{\rm ii)} The physical interpretation of  the optimal control Problem $\cQ^t$ is the following:
given a heat transfer process governed by the variational inequality (\ref{1d}) with the  data   $b$ and $q$ which satisfy condition (\ref{dw}) and (\ref{dwz}),
we are looking for an internal energy $f^*\in L^2(\Omega)$ such that
the temperature  $u$   is as close as possible to the ``desired" temperature $\phi$. Furthermore, this choice has to fulfill a minimum expenditure condition which is taken into account by the last  term in the cost functional. In fact, a compromise policy between the two aims (``$u$ close to $\phi$" and ``minimal energy $f$") has to be found and the relative importance of each criterion with respect to the other is expressed by the choice of the weight coefficients $\omega$ and  $\delta$.
Theorem \ref{t7} guarantees the  existence of at least one optimal energy function $f^*$  and, if the target $\phi$ vanishes, the optimal energy is unique. A similar comment can be made on the optimal control Problem $\cQ^t_n$.
Finally, the optimal solutions of the  heat transfer problem converge (in the sens given by Theorem \ref{t7} c)) to an optimal solution of the thermal problem with prescribed temperature on $\Gamma_2$, as the heat transfer coefficient converges to infinity.

\section*{Acknowledgement}

\indent This project has received funding from the European Union’s Horizon 2020
Research and Innovation Programme under the Marie Sklodowska-Curie
Grant Agreement No 823731 CONMECH.

\end{document}